\documentclass[11pt]{amsart}

\usepackage{color}

\usepackage{amsmath,amssymb,amsthm,amsfonts,amsbsy,amstext}

\usepackage{amsmath, amssymb} 
\usepackage{amsthm, amsfonts, mathrsfs}

\usepackage{amsfonts,graphicx}

\theoremstyle{plain}
\newtheorem{theo}{Theorem}[section]
\newtheorem{lemm}[theo]{Lemma}

\newtheorem{prop}[theo]{Proposition}

\theoremstyle{plain}
\theoremstyle{definition}

\theoremstyle{remark}
\newtheorem{rema}[theo]{Remark}
\newtheorem*{rema*}{Remark}

\newcommand{\ZZ}{\mathbb{Z}}  
\newcommand{\CC}{\mathbb{C}}  
\newcommand{\NN}{\mathbb{N}} 
\newcommand{\V}{\zeta}
\newcommand{\T}{\vartheta}
\newcommand{\RR}{\mathbb{R}} 
\newcommand{\EE}{\varepsilon} 
\newcommand{\DD}{\textnormal{D}} 
\newcommand{\XX}{X_{t,|\textnormal{D}|}} 
\newcommand{\PP}{\mathfrak{p}}

\newcommand{\QQ}{\mathbb{Q}} 
\newcommand{\HH}{\ell} 
\newcommand{\MM}{\mathfrak{m}}
\newcommand{\rr}{\mathcal{R}}
\numberwithin{equation}{section}

\def\lap{\Delta}
\date{}

\begin{document}

\title[Global well-posedness for  Boussinesq System]
{Global well-posedness for  Euler-Boussinesq system with critical dissipation}
\author[T. Hmidi] {T. Hmidi}
\address{IRMAR, Universit\'e de Rennes 1\\ Campus de
Beaulieu\\ 35~042 Rennes cedex\\ France}
\email{taoufik.hmidi@univ-rennes1.fr}
\author[S. Keraani]{S. Keraani}
\address{IRMAR, Universit\'e de Rennes 1\\ Campus de
Beaulieu\\ 35~042 Rennes cedex\\ France}
\email{sahbi.keraani@univ-rennes1.fr}
\author[F. Rousset]{F. Rousset}
\address{IRMAR, Universit\'e de Rennes 1\\ Campus de
Beaulieu\\ 35~042 Rennes cedex\\ France}
\email{frederic.rousset@univ-rennes1.fr}


\begin{abstract}
In this paper we study a fractional diffusion Boussinesq model which  couples 
  the incompressible Euler equation for the velocity and 
  a  transport equation with fractional diffusion
  for the temperature. We prove   global well-posedness results.

\end{abstract}

\maketitle

\section{Introduction} 
Boussinesq systems of the  type
\begin{equation*}
\label{bintro}
\left\{ 
\begin{array}{ll} 
\partial_{t}v+v\cdot\nabla v+\nabla p=\theta e_{2} + \nu \mathcal{D}_{v} v \\ 
\partial_{t}\theta+v\cdot\nabla\theta= \kappa \mathcal{D}_{\theta} \theta\\
\textnormal{div}\,v=0\\
v_{| t=0}=v^{0}, \quad \theta_{| t=0}=\theta^{0}  
\end{array} \right.
\end{equation*}
are simple models widely used in the modelling of oceanic and atmospheric motions.
These models also appear in many other physical problems,  we refer for instance  to  \cite{Brenier,Cons-Ber} for more details.
Here, we focus on the two-dimensional case, the space variable $x= (x_{1}, x_{2})$
 is in $\mathbb{R}^2$, the velocity field $v$    is given by  $v=(v^1,v^2)$ and the
  pressure $p$ and the temperature $\theta$ are scalar functions.
   The factor $\theta e_{2}$ in the velocity equation,  the vector $e_2$ being  given by $(0,1)$, models
    the effect of gravity on  the fluid motion. The operator $\mathcal{D}_{v} $ and 
    $\mathcal{D}_{\theta}$ whose form may vary are used to  take into account the possible
    effects of diffusion and dissipation in the fluid motion, thus  the constants 
    $\nu \geq 0$, $\kappa \geq 0$ can be seen as the inverse of Reynolds numbers.

   Mathematically, the simplest model to study is the fully viscous model  when  $\nu>0$, $\kappa>0$ and
    $\mathcal{D}_{v}= \Delta$,  $\mathcal{D}_{\theta}= \Delta.$ The properties of the system
     are very similar to the one of the two-dimensional Navier-Stokes equation
      and  similar global well-posedness results can be obtained.

  The most difficult  model   for the mathematical study is the inviscid  one, i.e. when $\nu= \kappa = 0$.
   A local existence result of smooth solution can be proven  as for  symmetric hyperbolic quasilinear systems,
    nevertheless,  it is not known if smooth solutions can develop singularities in finite time.
  Indeed, the temperature $\theta$ is the solution of a transport equation  and the vorticity $\omega=
   \mbox{curl } v = \partial_{1} v^2 - \partial_{2} v^1$ solves the equation
   \begin{equation}
   \label{vortintro}
   \partial_{t} \omega + v \cdot \nabla \omega =  \partial_{1} \theta.
   \end{equation}
   The main difficulty  is that to get an $L^\infty$ estimate on $\omega$ which is crucial to prove
    global existence of smooth solutions for  Euler type  equation, one needs to 
     estimate $\int_{0}^T ||\partial_{1} \theta ||_{L^\infty}$ and,  unfortunately, no {\it a priori }   estimate
      on $\partial_{1} \theta $ is known.

   In order to understand the coupling between  the two equations in Boussinesq type systems,
      there have been many recent works studying  Boussinesq systems with partial viscosity i.e. with
       a viscous term only in one equation.
      For $\kappa>0$, $\nu = 0$  and $\mathcal{D}_{\theta}= \Delta$,  the question of global existence is solved recently in a series of papers. In \cite{Cha}, Chae proved the global existence and uniqueness for initial data $(v^0,\theta^0)\in H^s\times H^s,$ with $s>2,$ see also  \cite{Hou}. This result was  recently extended   in \cite{hk} by the two first authors
   to initial data   $v^0\in B_{p,1}^{\frac{2}{p}+1}$ and $\theta^0\in B_{p,1}^{-1+\frac2p}\cap L^r, r>2$. More recently, the study of global existence of Yudovich solutions for this system has been done  in \cite{dp1}. We  also mention that in \cite{dp},  Danchin and Paicu were able to construct global  strong   solutions (still for $ \kappa>0$, $ \nu =0$) for a dissipative term of the form $ \mathcal{D}_{\theta} = \partial_{11}\theta$ instead of $\Delta\theta.$  Recently the first author and Zerguine \cite{H-Z} proved  the  global well-posedness for
        fractional diffusion $\mathcal{D}_\theta = - |\DD|^{ \alpha}$  for  $\alpha\in]1,2[$
         where the operator  $|\DD|^{ \alpha}$  is defined by 
         $$
\mathcal{F}(\vert \DD\vert^{\alpha}u)(\xi)=|\xi|^\alpha (\mathcal{F}u)(\xi).
$$
In these works, the global existence result relies on the fact 
 that the only smoothing effect due to the transport-(fractional) diffusion equation
 $$ \partial_{t} \theta + v \cdot \nabla \theta +  |\DD|^{\alpha} \theta = 0$$
  governing the temperature is sufficient to counterbalance the amplification of  the vorticity. However the case $\alpha=1$ is not reached by their method. The  main reason is that this case can be seen as critical in 
   the previous  approaches  in the sense that  the smoothing effect for  the temperature
    equation does not provide   the  
    $ L^1_{T}(L^\infty)$ bound for $\partial_{1} \theta $ which seems needed to control
     the amplification of  the $L^\infty$ norm of the vorticity. Note that  such an estimate is nevertheless 
      almost true  since $\partial_{1} \theta$ can be estimated in the space  $\tilde{L}^1_{T}(L^\infty)$
       which has the same scaling 
      (see below for the definition). The aim of this  paper is  the study of
       the well-posedness for  this case,
       i.e., we focus on  the system
     \begin{equation}
     \label{Bouss}
\left\{ 
\begin{array}{ll} 
\partial_{t}v+v\cdot\nabla v+\nabla p=\theta e_{2}  \\ 
\partial_{t}\theta+v\cdot\nabla\theta +  |\DD| \theta = 0 \\
\textnormal{div}\,v=0\\
v_{| t=0}=v^{0}, \quad \theta_{| t=0}=\theta^{0}.  
\end{array} \right.
     \end{equation}
  
   Note that  we have taken $\kappa =1$ which is legitimate since we study global well-posedness issues.
    Indeed, we can always change the coefficient $\kappa$
    into $1$ by a change of scale. 
    
   We also point out that at first sight,  the system \eqref{Bouss}  contains the mathematical  difficulties 
    of  the critical quasigeostrophic
    equation introduced in \cite{cmt}
    \begin{equation}
    \label{quasi} \partial_{t} \theta + v \cdot \nabla \theta +  |\DD| \theta = 0, \quad v = \nabla^{\perp}|\DD|^{-1} \theta
     \end{equation}
      which was much  studied recently. Indeed,  in \eqref{Bouss} the link between
       $v$ and $\theta$ is not given by the Riesz transform but by a dynamical equation, 
        the first equation of \eqref{Bouss}. Nevertheless, from this velocity  equation one gets  that $v$ has basically
         the regularity of $\theta$ as in the quasigeostrophic equation. 
        The global well-posedness for \eqref{quasi} was obtained recently by  Kiselev, 
  Nasarov and Volberg \cite{knv}.  We also refer to the work \cite{CV} by Caffarelli and Vasseur
   about the regularity of weak solutions. Other discussions  can be found in \cite{ab-Hm,ccw,cw}.
    
    \
    
    The main result of this paper is a  global well-posedness result for the system \eqref{Bouss} (see section \ref{preliminaries} for the  definitions  and the basic properties of Besov spaces).
 \begin{theo}\label{theo1}
   Let $p\in]2, \infty[$, $v^0\in {B}_{\infty, 1}^{1}\cap \dot{W}^{1,p}$ be a divergence-free vector  field of $\RR^2$  and 
$\theta^0\in  B_{\infty,1}^0\cap L^p$. Then there exists a unique global solution $(v, \theta)$ to the system 
\eqref{Bouss} with 
\begin{equation*}
v\in L^\infty_{\textnormal{loc}}\big(\RR_+;{B}_{\infty, 1}^{1}\cap \dot{W}^{1,p}),
\qquad \theta\in L_{\textnormal{loc}}^{\infty}\big(\RR_+; B_{\infty,1}^0\cap L^p\big)\cap \widetilde L_{\textnormal{loc}}^1(\RR_+; {B}_{p, \infty}^{1}).
\end{equation*}
\end{theo}
A few remarks are in order.

\begin{rema}  If we take $\theta=0$ then the system \eqref{Bouss} is reduced to the well-known 2D incompressible Euler system.  It is well known that this system is globally well-posed in  $H^s$ for  $s>2$. The main argument for globalization is the  BKM criterion  \cite{bkm} ensuring that the development of finite time singularities for Kato's solutions is related to the blowup  of the $L^\infty$ norm of the vorticity near the maximal time existence.
In \cite{vis} Vishik  has  extended the  global existence  of strong solutions result   to initial data lying on the spaces  $B_{p,1}^{1+2/p}.$  
Notice that  these spaces   have the same scaling as  Lipschitz functions (the space which is relevant for the hyperbolic theory)  and in this sense they are called  critical.  We emphasize that the application of the BKM criterion requires  a  super-lipschitzian regularity ($H^s$ with $s>2$ for example). For that reason the question of global existence in the critical spaces $B_{p,1}^{1+2/p}$ is  hard to deal  with because  these spaces have only a  lipschitzian regularity and the  BKM criterion cannot be used.   

Since 
$
 B_{p,1}^{1+2/p}\hookrightarrow {B}_{\infty, 1}^{1}\cap \dot{W}^{1,r}$ for all $ p\in [1,+\infty[$ and $r>\max\{p,2\}$, then 
  the space of initial velocity in our theorem  contains all the critical spaces $B_{p,1}^{1+2/p}$ except the biggest  one, that is  $B_{\infty,1}^{1}$. 
 For the limiting case  we have been  able to prove the global existence only  up to  the extra assumption   $\nabla v^0\in L^p$ for some $p\in]2, \infty[$.  The reason behind this extra assumption is the fact that to  obtain a global  $L^\infty$ bound for the vorticity we need before to establish an $L^p$ estimate for some $p\in]2,\infty[$ and it is not clear how to get rid of  this condition.
\end{rema}
\begin{rema} Since  $\nabla v, \nabla\theta\in L^{1}_{\rm loc}(\RR_+;L^\infty)$ (see Remark  \ref{rema7} below for $\theta$)  then we can easily propagate all the higher regularities: critical ({\it i.e.} $v_0\in B_{p,1}^{1+2/p}$ with $p$ finite) and sub-critical (for example  $v_0\in H^{s}$, for $s>2$).

\end{rema}

The main idea in the proof of Theorem \ref{theo1} is to really used the  structural properties of the system
 solved by $(\omega, \theta)$, $\omega = \mbox{curl } v= \partial_{1}v^2 - \partial_{2} v^1$.
  Indeed, if we neglect the nonlinear terms for the moment, one gets  the system
  $$ \partial_{t} \omega = \partial_{1} \theta, \quad \partial_{t} \theta = - |\DD| \theta$$
   and we notice that its symbol given by 
  $$ \mathcal{A}(\xi)  =   \left( \begin{array}{cc}   0 & i \xi_{1}  \\  0 & - |\xi| \end{array} \right) $$
   is diagonalizable  for $\xi \neq 0$ with two real distinct eigenvalues 
    which are  $0$ and $-|\xi|$.
      By using the Riesz transform $\mathcal{R}={ \partial_{1} / |\DD| } $,
       one gets  that the diagonal form of the system is given by
    $$ \partial_{t} \mathcal{R}  \theta = |\DD| \mathcal{R}  \theta, \quad
     \partial_{t}\big( \omega + \mathcal{R}\theta \big) = 0.$$
     This last form of the system is much more convenient in order to perform
       {\it a priori }   estimates.
    To prove Theorem \ref{theo1},  we shall use the same idea, we shall diagonalize the
     linear part of the system  and  then get  {\it a priori }   estimates from the study of  the new system.
   The main technical difficulty in  this program when one takes the nonlinear terms into account
    is to evaluate in a sufficiently sharp way the commutator
     $[\mathcal{R}, v \cdot \nabla ] $ between
     the  Riesz transform and  the convection operator. Such commutator estimates are  stated and proven 
      in section \ref{sectionRiesz} of the paper.
      
      The  diagonalization  approach used  in this paper also allows to prove global well-posedness
       in different spaces 
       for a "Boussinesq-Navier-Stokes" system i.e. the system which corresponds
        to $\mathcal{D}_{\theta}=0$ and $\mathcal{D}_{v} = -|\DD|$. This is discussed in a companion
        \mbox{ paper  \cite{HKR}.}
      
     \ 
          
   The  remaining of the paper is organized as follows.
  In section \ref{preliminaries} we recall some  functional spaces and
   we give some of their useful properties.
   Section \ref{sectionRiesz} is devoted to  the study of 
   some  commutators involving
    the Riesz transform. In section \ref{sectionTD}  we  study a linear transport-(fractional) diffusion equation. Especially, we  establish  some smoothing effects and a logarithmic estimate type.
    
    In section \ref{sectionbouss} we give the proof of Theorem \ref{theo1} which is 
     It is splitted into three parts.  We first establish some {\it a priori }   estimates, then we prove the uniqueness and finally  we briefly        explain how one can easily combine a procedure
      of smoothing out of the initial data with 
       the  {\it a priori }   estimates  to get the existence part of the theorem. An appendix  is devoted to the proof of  a technical              commutator lemma.

    \section{Notations and preliminaries}
    
    \label{preliminaries}
       \subsection{Notations}Throughout this work we will use the following notations.
       
$\bullet$ For any positive  $A$ and $B$  the notation  $A\lesssim B$ means that there exist a positive  harmless constant $C$ such that $A\le CB$. 

$\bullet$ For any tempered distribution $u$  both $ \hat u$  and $\mathcal F u$ denote the Fourier transform of $u$.

$\bullet$ Pour every $p\in [1,\infty]$, $\|\cdot\|_{L^p}$ denotes the norm in the Lebesgue space $L^p$.

$\bullet$ The norm  in the mixed space time Lebesgue space $L^p([0,T],L^r(\mathbb R^d)$ is denoted by  $\|\cdot\|_{L^p_TL^r}$ (with the obvious generalization to  $\|\cdot\|_{L^p_T\mathcal X} $ for any normed space $\mathcal X$).

$\bullet$ For any pair of operators $P$ and $Q$ on some Banach space $\mathcal{X}$, the commutator $[P,Q]$ is given by $PQ-QP$.

$\bullet$ For $p\in[1,\infty]$, we denote by $\dot{W}^{1,p}$ the space of distributions $u$ such that $\nabla u\in L^p.$

    \subsection{Functional spaces}

  Let us introduce the so-called   Littlewood-Paley decomposition and the corresponding  cut-off operators. 
There exists two radial positive  functions  $\chi\in \mathcal{D}(\RR^d)$ and  $\varphi\in\mathcal{D}(\RR^d\backslash{\{0\}})$ such that
\begin{itemize}
\item[\textnormal{i)}]
$\displaystyle{\chi(\xi)+\sum_{q\geq0}\varphi(2^{-q}\xi)=1}$;$\quad \displaystyle{\forall\,\,q\geq1,\, \textnormal{supp }\chi\cap \textnormal{supp }\varphi(2^{-q})=\varnothing}$\item[\textnormal{ii)}]
 $ \textnormal{supp }\varphi(2^{-j}\cdot)\cap
\textnormal{supp }\varphi(2^{-k}\cdot)=\varnothing,$ if  $|j-k|\geq 2$.
\end{itemize}

For every $v\in{\mathcal S}'(\RR^d)$ we set 
  $$
\Delta_{-1}v=\chi(\hbox{D})v~;\, \forall
 q\in\NN,\;\Delta_qv=\varphi(2^{-q}\hbox{D})v\quad\hbox{ and  }\;
 S_q=\sum_{ j=-1}^{ q-1}\Delta_{j}.$$
 The homogeneous operators are defined by
 $$
 \dot{\Delta}_{q}v=\varphi(2^{-q}\hbox{D})v,\quad \dot S_{q}v=\sum_{j\leq q-1}\dot\Delta_{j}v,\quad\forall q\in\ZZ.
 $$
From   \cite{b}  we split the product 
  $uv$ into three parts: $$
uv=T_u v+T_v u+R(u,v),
$$
with
$$T_u v=\sum_{q}S_{q-1}u\Delta_q v,\quad  R(u,v)=\sum_{q}\Delta_qu\tilde\Delta_{q}v  \quad\hbox{and}\quad \tilde\Delta_{q}=\sum_{i=-1}^1\Delta_{q+i}.
$$

\

 For
 $(p,r)\in[1,+\infty]^2$ and $s\in\RR$ we define  the inhomogeneous Besov \mbox{space $B_{p,r}^s$} as 
the set of tempered distributions $u$ such that
$$\|u\|_{B_{p,r}^s}:=\Big( 2^{qs}
\|\Delta_q u\|_{L^{p}}\Big)_{\ell ^{r}}<+\infty.
$$
The homogeneous Besov space $\dot B_{p,r}^s$ is defined as the set of  $u\in\mathcal{S}'(\RR^d)$ up to polynomials such that
$$
\|u\|_{\dot B_{p,r}^s}:=\Big( 2^{qs}
\|\dot\Delta_q u\|_{L^{p}}\Big)_{\ell ^{r}(\ZZ)}<+\infty.
$$
{Let $T>0$} \mbox{and $\rho\geq1,$} we denote by $L^\rho_{T}B_{p,r}^s$ the space of distributions $u$ such that 
$$
\|u\|_{L^\rho_{T}B_{p,r}^s}:= \Big\|\Big( 2^{qs}
\|\Delta_q u\|_{L^p}\Big)_{\ell ^{r}}\Big\|_{L^\rho_{T}}<+\infty.$$
We say that 
$u$ belongs to the space
 $\widetilde L^\rho_{T}{B_{p,r}^s}$ if
 $$
 \|u\|_{ \widetilde L^\rho_{T}{B_{p,r}^s}}:= \Big( 2^{qs}
\|\Delta_q u\|_{L^\rho_{T}L^p}\Big)_{\ell ^{r}}<+\infty.$$
By  a direct application  of 
the Minkowski inequality, we have the following links between these spaces. 
Let $ \varepsilon>0,$ then 
$$
L^\rho_{T}B_{p,r}^s\hookrightarrow\widetilde L^\rho_{T}{B_{p,r}^s}\hookrightarrow{L^\rho_{T}}{B_{p,r}^{s-\varepsilon}}
,\,\textnormal{if}\quad  r\geq \rho,$$
$$
{L^\rho_{T}}{B_{p,r}^{s+\varepsilon}}\hookrightarrow\widetilde L^\rho_{T}{B_{p,r}^s}\hookrightarrow L^\rho_{T}B_{p,r}^s,\, \textnormal{if}\quad 
\rho\geq r.
$$
 We will  make continuous use of Bernstein inequalities (see  \cite{che1} for instance).
\begin{lemm}\label{lb}\;
 There exists a constant $C$ such that for $q,k\in\NN,$ $1\leq a\leq b$ and for  $f\in L^a(\RR^d)$, 
\begin{eqnarray*}
\sup_{|\alpha|=k}\|\partial ^{\alpha}S_{q}f\|_{L^b}&\leq& C^k\,2^{q(k+d(\frac{1}{a}-\frac{1}{b}))}\|S_{q}f\|_{L^a},\\
\ C^{-k}2^
{qk}\|{\Delta}_{q}f\|_{L^a}&\leq&\sup_{|\alpha|=k}\|\partial ^{\alpha}{\Delta}_{q}f\|_{L^a}\leq C^k2^{qk}\|{\Delta}_{q}f\|_{L^a}.
\end{eqnarray*}

\end{lemm}

\section{Riesz transform and commutators}

\label{sectionRiesz}
 In the next proposition we gather some properties of  the Riez operator $\mathcal{R}={\partial_{1} }/{  |D|}.$
    \begin{prop}
    \label{cor1}  Let  $\mathcal{R}$ be the Riez operator $\mathcal{R}={\partial_{1} }/{  |D|}.$  Then the following hold true.
   \begin{enumerate}
 \item  For every $p \in ]1, +\infty[,$ 
    \begin{equation*}
  \label{CZ}
   \| \mathcal{R}  \|_{\mathcal{L}(L^p)} \lesssim 1.
   \end{equation*}
  \item Let $\chi\in\mathcal{D}(\RR^d)$. Then, there   exists $C>0$ such that 
$$
\||\textnormal{D}|^s\chi(2^{-q}|\textnormal{D}|)\mathcal{R}\|_{\mathcal{L}(L^p)}\le C 2^{qs},
$$
for every $(p,s,q)\in[1,\infty]\times ]0,+\infty[\times\NN$.

We notice that the previous results hold true if we change $|\textnormal{D}|^s$ by $\nabla^s$ \mbox{with $s\in\NN$.}\\
\item Let $\mathcal C$ be a fixed ring. Then, there exists $\psi\in \mathcal S$
whose spectum does not meet the origin such that 
$$
\mathcal{R} f=2^{qd}\psi(2^{q}\cdot)\star f
$$
for every $f$ with Fourier transform supported  in $2^q\mathcal C$.
\end{enumerate}
\end{prop}   
\begin{proof}
(1) It is a classical Calder\'on-Zygmund theorem (see \cite{Stein} for instance).

(2) Let $K\in\mathcal{S}'$  such that  $\mathcal{F} K(\xi)=|\xi|^s\chi(|\xi|)\frac{\xi_i}{|\xi|}.$   $K$ is a  tempered distribution such that its Fourier transform  $\mathcal{F} K$ is  \mbox{$\mathcal{C}^{\infty}(\mathbb{R}^d\setminus\{0\})$}  and satisfies for every $\alpha\in\NN^d$
 $$
 |\partial _{\xi}^{\alpha}\mathcal{F}K(\xi)|\leq C_{\alpha}|\xi|^{s-|\alpha|}.
 $$ 
According to Mikhlin-H\"ormander  Theorem  (see \cite{Stein} for instance)
we have 
  for every  $\alpha\in\NN^d$ and $x\in\RR^d\setminus\{0\}$, 
$$|\partial_x^{\alpha}K(x)|\leq C_{\alpha}'|x|^{-s -d-|\alpha|}.$$
 Then, it ensues that 
$$
|K(x)|\leq C|x|^{-d-s}\qquad \forall x\neq 0.
$$
Since  $\mathcal{F}K$ is obviously in  $L^1$, we also have that  $K\in \mathcal{C}_0(\RR^d)$. This removes the singularity at the origin and gives  
$$
|K(x)|\le C(1+|x|)^{-d-s}\qquad \forall x\in\mathbb R^d.
$$ 
Therefore we get that the kernel $K\in L^1$. Now for every $u\in L^p$ we have $|\textnormal{D}|^s\mathcal{R}\chi(|\textnormal{D}|) u=K\star u.$ We   can now  conclude the case $q=0$  by  using the  classical
 Young inequality for  convolution products.\\
The case $q\geq 1$  can be derived from  $q=0$ via an obvious argument of  homogeneity.

(3) This is can be done easily by introducing  a judiciously chosen  cut-off function.
\end{proof}

The following  lemma  will be useful in the proof of many  commutator estimates.
 \begin{lemm}
 \label{commu}Given $(p,m)\in[1,\infty]^2$ such that $p\geq m'$ with $m'$ the conjugate exponent of $m$. Let $f,g$ and $h$ be three functions such that $\nabla f\in L^p, g\in L^m$ and $xh\in L^{m'}$. Then,
 $$
 \|h\star(fg)-f(h\star g)\|_{L^p}\leq \|xh\|_{L^{m'}}\|\nabla f\|_{L^p}\|g\|_{L^{m}}.
 $$

 \end{lemm}
 \begin{proof}
We have by definition and Taylor formula
\begin{eqnarray*}
h\star(fg)(x)-f(h\star g)(x)&=&\int_{\RR^d} h(x-y)g(y)\big(f(y)-f(x)  \big)dy\\
&=&\int_0^1\int_{\RR^d} g(y)h(x-y)\Big[(y-x)\cdot\nabla f(x+t(y-x))\Big]dydt.
\end{eqnarray*}
Using H\"{o}lder inequalities and making a change of variables $z=t(x-y)$ we get
$$
|h\star(fg)(x)-f(h\star g)(x)|\le \|g\|_{L^m}\int_0^1\Big(\int_{\RR^d} t^{-d}h_1^{m'}(t^{-1}z)|\nabla f|^{m'}(x-z)dz\Big)^{\frac{1}{m'}},
$$
where we set $h_1(z)=|z||h(z)|.$ Using Convolution inequalities we obtain since $p\geq m'$
\begin{eqnarray*}
\|h\star(fg)-f(h\star g)\|_{L^p}&\le& \|g\|_{L^m}\|h_1^{m'}\|_{L^1}^{\frac{1}{m'}}\||\nabla f|^{m'}\|_{L^{\frac{p}{m'}}}^{\frac{1}{m'}}\\
&\leq&\|g\|_{L^m}\|h_1\|_{L^{m'}}\|\nabla f\|_{L^p}.
\end{eqnarray*}

\end{proof}

As explained in the introduction,  the control of the commutator between $\mathcal{R}$ and the convection
 operator $v \cdot \nabla $ is a  crucial ingredient in the proof of Theorem \ref{theo1}.
\begin{theo}\label{propcom}
Let  $v$ be is  a smooth  divergence-free  vector field. 
\begin{enumerate}
\item For every $(p,r)\in [2,\infty[\times[1,\infty]$ there exists 
a constant $C=C(p,r)$ such that 
$$
\|[\mathcal{R}, v\cdot\nabla]\theta\|_{B_{p,r}^0}\le C \|\nabla v\|_{L^p}\big(\|\theta\|_{B_{\infty,r}^0}+\|\theta\|_{L^p}\big),
$$
for every smooth scalar function $\theta$.
\item 
For every  $(r,\rho)\in[1,\infty]\times]1,\infty[$  and $\epsilon >0$ there exists 
a constant $C=C(r,\rho,\varepsilon)$ such that 

 $$
 \|[\mathcal{R}, v\cdot\nabla]\theta\|_{B_{\infty,r}^0}\le C (\|\omega\|_{L^\infty}+\|\omega\|_{L^\rho})\big(\|\theta\|_{B_{\infty,r}^\epsilon}+\|\theta\|_{L^\rho}\big),
 $$
 for every smooth scalar function $\theta$.
 \end{enumerate}
 \end{theo}
 \begin{proof}
 We  split the commutator into three parts, according to  Bony's decomposition,
\begin{eqnarray*}
\nonumber[\mathcal{R}, v\cdot\nabla]\theta&=&\sum_{q\in\NN}[\mathcal{R}, S_{q-1}v\cdot\nabla]\Delta_q\theta+\sum_{q\in\NN}[\mathcal{R}, \Delta_qv\cdot\nabla]S_{q-1}\theta\\
\nonumber&+&\sum_{q\geq-1} [\mathcal{R}, \Delta_qv\cdot\nabla]\widetilde{\Delta}_q\theta\\
\nonumber&=& \sum_{q\in\NN}\mbox{I}_q+\sum_{q\in\NN}\mbox{II}_q+\sum_{q\geq-1}\mbox{III}_q\\
&=&\mbox{I}+\mbox{II}+\mbox{III}.
\end{eqnarray*}
We start with the estimate of the first term $\mbox{I}$.    According to the point (3) of \mbox{Proposition \ref{cor1}} there  exists  $h\in\mathcal{S}$  whose spectum does not meet the origin such that 
 $$
\mbox{I}_q(x)=h_q\star( S_{q-1}v\cdot\nabla\Delta_q\theta)-S_{q-1}v\cdot(h_q\star \nabla\Delta_q\theta),
$$
 where $h_q(x)=2^{dq}h(2^qx)$.
 Applying Lemma \ref{commu} with $m=\infty$   we get 
  \begin{eqnarray}
  \nonumber
\|\mbox{I}_q\|_{L^p}&\lesssim &\|xh_q\|_{L^1}  \|\nabla S_{q-1} v\|_{L^p}\|\Delta_q\nabla\theta\|_{L^\infty}
\\
\label{x1}
&\lesssim& \|\nabla v\|_{L^p}\|\Delta_q\theta\|_{L^\infty}.
\end{eqnarray}
In the last line we've used Bernstein inequality  and  $\|xh_q\|_{L^1}= 2^{-q}\|xh\|_{L^1}$.

 Combined with the trivial fact 
 $$\Delta_j\sum_{q}\mbox{I}_q= \sum_{|j-q|\le 4}\mbox{I}_q
 $$ this yields
\begin{eqnarray*}
\|\mbox{I}\|_{B_{p,r}^0}
&\lesssim&
\Big(\sum_{q\geq-1}\|\mbox{I}_q\|_{L^p}^r\Big)^{\frac1r}
\\
&\lesssim&\|\nabla v\|_{L^p}
\|\theta\|_{B_{\infty,r}^0}.
\end{eqnarray*}

Let us move to  the second term $\mbox{II}$. As before one writes
$$
\mbox{II}_q(x)=h_q\star( \Delta_q v\cdot\nabla  S_{q-1}\theta)-\Delta_q v\cdot(h_q\star \nabla S_{q-1}\theta),
$$ 
and then we obtain the estimate 
\begin{eqnarray*}
\nonumber\|\mbox{II}_q\|_{L^p}&\lesssim& 2^{-q}\|\Delta_q \nabla v\|_{L^p}\| S_{q-1}\nabla\theta\|_{L^\infty}\\
&\lesssim&\|\nabla v\|_{L^p}\sum_{j\le q-2}2^{j-q}\|\Delta_j\theta\|_{L^\infty}.
\end{eqnarray*}
Combined with convolution inequalities this  yields
$$
\|\mbox{II}\|_{B_{p,r}^0}\lesssim \|\nabla v\|_{L^p}\|\theta\|_{B_{\infty,r}^0}.
$$

\

Let us now deal with the third term  $\mbox{III}$. Using  that the divergence of $\Delta_q v$ vanishes,  we rewrite  $\mbox{III}$  as
\begin{eqnarray*}
\mbox{III}&=&\sum_{q\geq 2} \mathcal{R} \textnormal{div} (\Delta_qv\, \widetilde{\Delta}_q\theta)- \sum_{q\geq 2} \textnormal{div} (\Delta_qv\,\mathcal{R}\widetilde{\Delta}_q\theta)+\sum_{q\leq 1}[\mathcal{R}, \Delta_{q}v\cdot\nabla]\widetilde{\Delta}_{q}\theta\\&=&J_1+J_2+J_3.
\end{eqnarray*}
Using   Proposition \ref{cor1}-(2), we get 
\begin{eqnarray*}
\big\|\Delta_j\mathcal{R}\textnormal{div} (\Delta_qv\, \widetilde{\Delta}_q\theta)\big\|_{L^p}\lesssim  2^j\|\Delta_q v\|_{L^p}\|\widetilde{\Delta}_q\theta\|_{L^\infty}.
\end{eqnarray*}
Also, since $\widetilde{\Delta}_q\theta$ is supported away from zero for $q\geq 2$ then  Proposition \ref{cor1} (3) yields
\begin{eqnarray*}
 \big\|\Delta_j\textnormal{div} (\Delta_qv\, \mathcal{R}\widetilde{\Delta}_q\theta)\big\|_{L^p}&\lesssim & 2^j\|\Delta_q v\|_{L^p}\|\mathcal R\widetilde{\Delta}_q\theta\|_{L^\infty}
 \\
 &\lesssim& 2^j\|\Delta_q v\|_{L^p}\|\widetilde{\Delta}_q\theta\|_{L^\infty}.
\end{eqnarray*}
Therefore we get
\begin{eqnarray*}
\|\Delta_j (J_1+J_2)\|_{L^p}&\lesssim&  \sum_{q\in\NN\atop q\geq j-4}2^j\|\Delta_q v\|_{L^p}\|\widetilde{\Delta}_q\theta\|_{L^\infty}\\
&\lesssim&\|\nabla v\|_{L^p} \sum_{q\in\NN\atop q\geq j-4}2^{j-q}\|{\Delta}_q\theta\|_{L^\infty},
\end{eqnarray*}
where we have again used Bernstein inequality to get the last line.
It suffices now to use convolution inequalities to get
$$
\|J_1+J_2\|_{B_{p,r}^0}\lesssim \|\nabla v\|_{L^p}\|\theta\|_{B_{\infty,r}^0}.
$$
For the last term $J_3$ we can  write
$$
\sum_{-1\leq q\leq 1}[\mathcal{R}, \Delta_{q}v\cdot\nabla]\widetilde{\Delta}_{q}\theta(x)=\sum_{q\leq 1}[\textnormal{div }\widetilde{\chi}(\textnormal{D})\mathcal{R}, \Delta_{q}v]\widetilde{\Delta}_{q}\theta(x),$$
where  $\widetilde{\chi}$ belongs to $\mathcal{D}(\RR^d)$.
Proposition  \ref{cor1} ensures  that 
$\textnormal{div }\widetilde{\chi}(\textnormal{D})\mathcal{R}$ is a convolution operator  with  a kernel   $\tilde{h}$ satisfying 
$$
|\tilde{h}(x)|\lesssim (1+|x|)^{-d-1}.
$$  
Thus
$$
J_3=  \sum_{q\leq 1} \tilde h\star( \Delta_qv\cdot\tilde\Delta_q\theta)-\Delta_qv\cdot(\tilde h\star \tilde\Delta_q\theta).
$$
First of all we point out that $\Delta_j J_3=0$ for $j\geq 6$, thus we just need  to estimate the low frequencies of $J_3$.   Noticing that $x\tilde h$ belongs to $L^{p'}$ for  $p'>1$ then using Lemma \ref{commu} with $m=p\geq 2$ we obtain
\begin{eqnarray*}
\|\Delta_j J_3\|_{L^\infty}&\lesssim&  \sum_{q\le 1}  \|x\tilde h\|_{L^{p'}} \|\Delta_q \nabla v\|_{L^p}\|\widetilde{\Delta}_q\theta\|_{L^p}\\
&\lesssim& \| \nabla v\|_{L^p}\sum_{-1\leq q\leq 1}\|{\Delta}_q\theta\|_{L^p}.
\end{eqnarray*}
This yields  finally
$$
\|J_3\|_{B_{p,r}^0}\lesssim \|\nabla v\|_{L^p}\|\theta\|_{L^p}.
$$
This completes the proof of  the first part of Theorem \ref{propcom}. The second part can be done in the same way so that we will only give here a shorten  proof. To estimate the terms ${\rm I}$ and ${\rm II}$ we use two facts: the first one is $\|\Delta_q\nabla u\|_{L^\infty}\approx \|\Delta_q \omega\|_{L^\infty}$ for all  $q\in \mathbb N.$ The second one is 
 \begin{eqnarray*}
\|\nabla S_{q-1} v\|_{L^\infty}& \lesssim&\|\nabla\Delta_{-1}v\|_{L^\infty}+\sum_{j=0}^{q-2}\|\Delta_j\nabla v\|_{L^\infty}\\
&\lesssim&\|\omega\|_{L^\rho}+ q\|\omega\|_{L^\infty}.
\end{eqnarray*}
For the remainder term we do strictly the same analysis as before except for $J_3$: we apply Lemma \ref{commu} with $p=\infty$ and $m=\rho$ leading to 
\begin{eqnarray*}
\nonumber
\|\Delta_j J_3\|_{L^p}&\lesssim&  \sum_{q\le 1}  \|x\tilde h\|_{L^{\rho'}} \|\Delta_q \nabla v\|_{L^\infty}\|\widetilde{\Delta}_q\theta\|_{L^\rho}\\
&\lesssim& \| \nabla v\|_{L^\rho}\sum_{-1\leq q\leq 1}\|{\Delta}_q\theta\|_{L^\rho}\\
&\lesssim& \|\omega\|_{L^\rho}\|\theta\|_{L^\rho}.
\end{eqnarray*}
This ends  the proof of the theorem.

\end{proof}
\section{Transport-Diffusion equation}
\label{sectionTD}
In this section we will give some  useful estimates for any smooth  solution of  a linear  transport-diffusion  model   given by  

\begin{equation} 
\left\{ \begin{array}{ll} 
\partial_{t} \theta+v\cdot\nabla \theta+\vert\textnormal{D}\vert \theta=f\\
\theta_{| t=0}=\theta^{0}.
\end{array} \right. \tag{${\textnormal{TD}}$}
\end{equation} 
We will discuss  three kinds of estimates:  $L^p$ estimates,  smoothing effects and logarithmic estimates.  
 
The proof of the following  $L^p$ estimates can be found in \cite{cc}.
\begin{prop}
\label{propmax}
Let  $v$ be  a smooth  divergence-free vector field of $\RR^d$  and  $\theta$ be a smooth solution of $({\textnormal{TD}})$. Then we have for every $p\in[1,\infty]$
$$
\|\theta(t)\|_{L^p}\le\|\theta^0\|_{L^p}+\int_0^t\|f(\tau)\|_{L^p}d\tau.
$$
\end{prop}
We intend to prove the following smoothing effect. 
\begin{theo}
 \label{thm99} Let  $v$ be a smooth divergence-free vector field of $\RR^d$ with vorticity $\omega$. Then, for every   $p\in [1,\infty[$  there exists a constant $C$ such that 
$$
\sup_{q\in\NN}2^q\|\Delta_q \theta\|_{L^1_tL^p}\leq C\| \theta^0\|_{L^p}+C\|\theta^0\|_{L^\infty} \|\omega\|_{L^1_tL^p},
$$
for every  smooth  solution $\theta$ of ${\rm (TD)}$ with $f\equiv 0$.
\end{theo}
\begin{proof}
We start with localizing in frequencies the equation:
for $q\geq-1$ we set $\theta_q:=\Delta_q\theta. $ Then 
$$
\partial_t\theta_q+v\cdot\nabla\theta_q+\vert\textnormal{D}\vert\theta_q=-[\Delta_q, v\cdot\nabla]\theta.
$$
Recall that ${\theta}_q$ is real function since the functions involved in the dyadic partition of the unity are radial. Then multiplying the above equation  by $|\theta_q|^{p-2}{\theta}_q,$ integrating by parts  and using H\"{o}lder inequalities we get
$$
\frac1p\frac{d}{dt}\|\theta_q\|_{L^p}^p+\int_{\RR^2}(\vert\textnormal{D}\vert\theta_q) |\theta_q|^{p-2}{\theta}_qdx\leq \|\theta_q\|_{L^p}^{p-1}\|[\Delta_q, v\cdot\nabla]\theta\|_{L^p}.
$$
Recall from \cite{cmz} the following generalized Bernstein inequality
$$
c  2^{q}\|\theta_q\|_{L^p}^p\le\int_{\RR^2}(\vert\textnormal{D}\vert\theta_q) |\theta_q|^{p-2}{\theta}_qdx,
$$
where $c$ depends on $p$. Inserting this estimate in the previous one  we obtain 
$$
\frac1p\frac{d}{dt}\|\theta_q\|_{L^p}^p+c2^q \|\theta_q\|_{L^p}^p\lesssim  \|\theta_q\|_{L^p}^{p-1} \|[\Delta_q, v\cdot\nabla]\theta\|_{L^p}.
$$
Thus we find 
\begin{equation}
\label{est}
\frac{d}{dt}\|\theta_q\|_{L^p}+c2^q \|\theta_q\|_{L^p}\lesssim  \|[\Delta_q, v\cdot\nabla]\theta\|_{L^p}.
\end{equation}
To estimate the right hand-side, we shall use the following lemma (see the appendix for the proof of this lemma).
\begin{lemm}
\label{propDelta}
Let $v$ be a smooth divergence-free vector field and $\theta$ be a smooth scalar function. Then,  for all $ p\in[1,\infty]$ and $ q\geq -1,$
$$
\|[\Delta_q, v\cdot\nabla]\theta\|_{L^p}\lesssim \|\nabla v\|_{L^p}\|\theta\|_{B_{\infty,\infty}^{0}}.
$$
\end{lemm}
Combined with \eqref{est} this lemma yields
\begin{eqnarray*}
\frac{d}{dt}\big( e^{ct2^q}\|\theta_q(t)\|_{L^p}\big)&\lesssim&   e^{ct2^q} \|\nabla v(t)\|_{L^p}\|\theta(t)\|_{B_{\infty,\infty}^0}\\
&\lesssim& e^{ct2^q}  \|\omega(t)\|_{L^p}\|\theta^0\|_{L^\infty}.
\end{eqnarray*}
To get the last line, we have used   the conservation of the $L^\infty$ norm of $\theta$ and the classical fact 
$$
 \|\nabla v \|_{L^p } \lesssim \| \omega \|_{L^p}\qquad\forall p\! \in ]1, +\infty[.
$$
Integrating the differential inequality  we get
\begin{eqnarray*}
\|\theta_q(t)\|_{L^p}&\lesssim& \|\theta_q^0\|_{L^p} e^{-ct2^q}+\|\theta^0\|_{L^\infty}\int_0^t e^{-c(t-\tau)2^q} \|\omega(\tau)\|_{L^p}d\tau.
\end{eqnarray*}
Integrating in time yields finally
\begin{eqnarray*}
2^q\|\theta_q\|_{L^1_tL^p}&\lesssim& \|\theta_q^0\|_{L^p} +\|\theta^0\|_{L^\infty}\int_0^t  \|\omega(\tau)\|_{L^p}d\tau\\
&\lesssim&\|\theta^0\|_{L^p} +\|\theta^0\|_{L^\infty}\int_0^t  \|\omega(\tau)\|_{L^p}d\tau,
\end{eqnarray*}
which is the desired result.
\end{proof}
Let us now move to the last part of this section which deals with some logarithmic estimates generalizing  the results of \cite{vis,hk}. First we recall 
the following  result of propagation of Besov regularities.
\begin{prop}\label{prop-Bes}
Let $(p,r)\in[1,\infty]^2, s\in]-1,1[$ and $\theta$ a smooth solution of {\rm ({\rm TD})}. Then we have
$$
\|\theta\|_{\widetilde L^\infty_tB_{p,r}^s}\lesssim e^{CV(t)}\Big(\|\theta^0\|_{B_{p,r}^s}+\int_0^te^{-CV(\tau)}\|f(\tau)\|_{B_{p,r}^s}d\tau \Big),
$$
where $V(t)=\|\nabla v\|_{L^1_t L^\infty}.$
\end{prop}
The proof of this result is omitted here and it can be done similarly to the inviscid case \cite{che1},   using especially  Proposition \ref{propmax}.

Now we will  show  that for the index regularity $s=0$ we can obtain a better estimate with a linear growth on Lipschitz norm of the velocity. 
\begin{theo}
\label{thmlog}
There exists $C>0$ such that if $\kappa\geq 0$, $ p\in[1,\infty]$ and  $\theta$  a solution of
$$
(\partial_t+v\cdot\nabla+\kappa|\DD|)\theta=f,
$$
then we have
$$
\|\theta\|_{\widetilde L^\infty_tB_{p,1}^0}\leq C\Big( \|\theta^0\|_{B_{p,1}^0}+\|f\|_{L^1_tB_{p,1}^0}\Big)\Big( 1+\int_0^t\|\nabla v(\tau)\|_{L^\infty}d\tau \Big).
$$
\end{theo}
\begin{proof}
We mention that the result is first proved in \cite{vis} for the case $\kappa=0$ by using the special structure of the transport equation. In \cite{H-K2} the first two authors  generalized Vishik's result for a transport-diffusion equation where the dissipation term has the form $-\kappa\Delta\theta$.  The method described in \cite{H-K2} can be easily adapted  here for  our model. 

 Let 
$q\in\NN\cup\{-1\}$ and  denote by  $\overline\theta_q$  the unique global solution of the initial value problem
\begin{equation}\label{R_{Q}}\left\lbrace
\begin{array}{l}
\partial_t \overline\theta_q+v\cdot\nabla \overline\theta_q+|\DD|\overline\theta_q=\Delta_{q}f,
\\
{\overline\theta_q}_{|t=0}=\Delta_{q}\theta^{0}.\\
\end{array}
\right.
\end{equation}
Using  Proposition \ref{prop-Bes}  with $s=\pm\frac12$   we get
$$
\|\overline\theta_q\|_{\widetilde L^\infty_tB_{p,\infty}^{\pm\frac12}}\lesssim
\big(\|\Delta_{q} \theta^0\|_{B_{p,\infty}^{\pm\frac12}}+
\|\Delta_{q}f\|_{L^1_{t}B_{p,\infty}^{\pm\frac12}}\big) e^{CV(t)
},
$$
where $
V(t)=\|\nabla v\|_{L^1_{t}L^\infty}.
$
Combined with  the  definition of Besov spaces this yields \mbox{ for $j,q\geq-1$}
\begin{equation}
\label{t1}
\|\Delta_{j}\overline\theta_q\|_{L^\infty_tL^p}\lesssim 2^{-\frac12|j-q|}
\big(\|   \Delta_{q} \theta^0\|_{L^p}+\|\Delta_{q}f\|_{L^1_{t}L^p}\big) e^{CV(t)}.
\end{equation}
By linearity and again  the definition of Besov spaces we have
\begin{eqnarray}\label{t2}
\|\theta\|_{\widetilde L^\infty_tB_{p,1}^0}\leq
\sum_{|j-q|\geq N}
\|\Delta_{j}\overline\theta_q\|_{L^\infty_tL^p}+\sum_{|j-q|< N}
\|\Delta_{j}\overline\theta_q\|_{L^\infty_tL^p},
\end{eqnarray}
where $N\in \Bbb N$ is to be chosen later.
To deal with the first sum we use (\ref{t1}) 
\begin{eqnarray*}
\nonumber\sum_{|j-q|\geq N}
\|\Delta_{j}\overline\theta_q\|_{L^\infty_tL^p}&\lesssim& 2^{-N/2}\sum_{q\geq-1}\big(\|\Delta_{q}\theta^0\|_{L^p}+\|\Delta_{q}f\|_{L^1_{t}L^p}\big)e^{CV(t)}
\\
&\lesssim&  2^{-N/2}
\big(\|\theta ^0\|_{B_{p,1}^0}+\|f\|_{L^1_{t}B_{p,1}^0}\big)e^{CV(t)}.
\end{eqnarray*}
We now turn to the second  sum in  the right-hand side of (\ref{t2}).  

It is clear that
\begin{equation*}\nonumber\sum_{|j-q|< N}
\|\Delta_{j}\overline\theta_q\|_{L^\infty_tL^p}\lesssim \sum_{|j-q|< N}
\|\overline\theta_q\|_{L^\infty_tL^p}.
\end{equation*}
Applying Proposition \ref{propmax} to  the system (\ref{R_{Q}}) yields 
$$
\|\overline\theta_q\|_{L^\infty_tL^p}\leq \|\Delta_q\theta^0\|_{L^p}+\|\Delta_{q}f\|_{L^1_{t}L^p}.
$$
It follows that 
\begin{equation*}
\sum_{|j-q|<N}\|\Delta_{j}\overline\theta_q\|_{L^\infty_tL^p}\lesssim N\big(\|\theta ^0\|_{B_{p,1}^0}+\|f\|_{L^1_{t}B_{p,1}^0}\big).
\end{equation*}
The outcome is  the following
$$
\|\theta\|_{\widetilde L^\infty_tB_{p,1}^0}\lesssim\big(\|\theta ^0\|_{B_{p,1}^0}+\|f\|_{L^1_{t}B_{p,1}^0}\big)\Big(2^{- N/2}e^{CV(t)}+N\Big).
$$
Choosing
$$
N=\Big[\frac{2C
V(t)}{ \log 2}\Big]+1,
$$
we get the desired result.

\end{proof}


\section{Proof of Theorem \ref{theo1} }

Throughout this section we use the notation $\Phi_k$ to denote any function
of the form 
$$
\Phi_k(t)=  C_{0}\underbrace{ \exp(...\exp  }_{k\,times}(C_0t)...),
$$
where $C_{0}$ depends on the involved norms of the initial data and its value may vary from line to line up to some absolute constants. 
We will make an intensive  use (without mentionning it) of  the following trivial facts
$$
\int_0^t\Phi_k(\tau)d\tau\leq \Phi_k(t)\qquad{\rm and}\qquad \exp({\int_0^t\Phi_k(\tau)d\tau})\leq \Phi_{k+1}(t).
$$

\label{sectionbouss}
The proof of Theorem \ref{theo1} will be done in several steps. The first one deals with  some {\it {\it a priori }   } estimates for the equations  \eqref{Bouss}. In the  second one we prove the uniqueness part. Finally,   we will  discuss    the construction of the solutions at the end of this section.
\subsection{\it a priori    estimates}
As we will see the important quantities to bound for all times are $\|\omega(t)\|_{L^\infty}$ and $\|\nabla v(t)\|_{L^\infty}$. It seems that for subcritical regularities like for example $ H^s, s>2$ or more generally  $B_{p,r}^{s},\, s>1+\frac2p$ we need only to bound the quantity $\|\partial_1\theta\|_{L^1_t L^\infty}$, which in turn controls  $\|\omega(t)\|_{L^\infty},$ due to Brezis-Gallo\"{u}et logarithmic estimate, see for example \cite{Shu}.  Even though these quantities seem to be less regular than $\|\nabla v\|_{L^1_t L^\infty}$, it is not at all clear how to estimate them without involving the latter quantity.  

When we deal with critical regularities which is our subject here  one needs to bound the Lipschitz norm of the velocity and this will require some refinement analysis, especially Theorem  \ref{thmlog} seems to be very crucial.  To obtain  a Lipschitz bound  we will proceed in several steps: one of the main step is to  give an $L^\infty$-bound of the vorticity but due to some technical difficulties related to Riesz transforms  this will be not done in a straight way. We prove before  an $L^p$ estimate  with $2<p<\infty$  and  this allows us to bound the vorticity in $L^\infty$.    

 We start with recalling the  $L^p$ estimate for the temperature function. It is a direct consequence of Proposition \ref{propmax}. 
 \begin{prop}
 \label{prop5.1}
Let $(v,\theta)$ be a smooth solution of \eqref{Bouss}, then for all $p\in [1,+\infty]$
$$
\|\theta(t)\|_{L^p}\leq \|\theta^0\|_{L^p}.
$$
\end{prop}
We intend  now  to bound  the $L^p$-norm of the vorticity and to describe a smoothing effect for the temperature.
\begin{prop}\label{max-pro}
If $\omega^0\in L^p$ and  $\theta^0\in L^p\cap L^\infty$ with $p\in]2,\infty[,$  then
$$
\|\omega(t)\|_{L^p}+\|\theta\|_{\widetilde L^1_t B_{p,\infty}^1}\le \Phi_1(t).
$$
\end{prop}
\begin{proof}
 Applying Riesz transform $\mathcal{R}$ to the temperature equation  we get
\begin{equation}\label{bm}
\partial_t\mathcal{R}\theta+v\cdot\nabla\mathcal{R}\theta+|\textnormal{D}|\mathcal{R}\theta=-[\mathcal{R}, v\cdot\nabla]\theta.
\end{equation}
Since $|\textnormal{D}|\mathcal{R}=\partial_1$ then the function $\Gamma:=\omega+\mathcal{R}\theta$ satisfies 
\begin{equation}
\label{bm1}
\partial_t\Gamma+v\cdot\nabla\Gamma=-[\mathcal{R}, v\cdot\nabla]\theta.
\end{equation}
According to  the first part of Theorem  \ref{propcom} applied with    $r=2$ we have
$$
\big\|[\mathcal{R},v\cdot\nabla]\theta\big\|_{B_{p,2}^0}\lesssim\|\nabla v\|_{L^p}\big(\|\theta\|_{B_{\infty,2}^0}+\|\theta\|_{L^p} \big).
$$
Using the classical   embedding $B_{p,2}^0\hookrightarrow L^p$ which is  true only for $p\in[2,\infty)$  
$$
\big\|[\mathcal{R},v\cdot\nabla]\theta\big\|_{L^p}\leq\|\nabla v\|_{L^p}\big(\|\theta\|_{B_{\infty,2}^0}+\|\theta\|_{L^p}\big).
$$
Since ${\rm div }\,v=0$ then we get from  the transport equation \eqref{bm1}
$$
\|\Gamma(t)\|_{L^p}\leq\|\Gamma^0\|_{L^p}+\int_0^t\|[\mathcal{R}, v\cdot\nabla]\theta(\tau)\|_{L^p}d\tau.
$$
Putting together the last two estimates  we get
\begin{eqnarray*}
\|\Gamma(t)\|_{L^p}&\lesssim& \|\Gamma^0\|_{L^p}+\int_0^t\|\nabla v(\tau)\|_{L^p}\big(\|\theta(\tau)\|_{B_{\infty,2}^0}+\|\theta\|_{L^p}\big) d\tau
\\
&\lesssim&\|\omega^0\|_{L^p}+\|\theta^0\|_{L^p}+\int_0^t\|\omega(\tau)\|_{L^p}\big(\|\theta(\tau)\|_{B_{\infty,2}^0}+\|\theta^0\|_{L^p}\big) d\tau.
\end{eqnarray*}
We have used here the Calder\'on-Zygmund estimates: for $p\in(1,\infty)$
$$
\|\nabla v\|_{L^p}\le C\|\omega\|_{L^p}\quad\hbox{and} \quad\|\mathcal{R}\theta^0\|_{L^p}\leq C\|\theta^0\|_{L^p}.
$$
On the other hand, from the continuity of the Riesz transform and   Proposition \ref{prop5.1}
\begin{eqnarray*}
\|\omega(t)\|_{L^p}&\le&\|\Gamma(t)\|_{L^p}+\|\mathcal{R}\theta\|_{L^p}\\
&\lesssim&\|\Gamma(t)\|_{L^p}+\|\theta^0\|_{L^p}.
\end{eqnarray*}
This leads to 
$$
\|\omega(t)\|_{L^p}\lesssim \|\omega^0\|_{L^p}+\|\theta^0\|_{L^p}+\int_0^t\|\omega(\tau)\|_{L^p}\big(\|\theta(\tau)\|_{B_{\infty,2}^0}+\|\theta^0\|_{L^p}\big)d\tau.
$$
According to Gronwall lemma we get
\begin{equation}
\label{gr1}
\|\omega(t)\|_{L^p}\le C_0e^{C_0 t}e^{C\|\theta\|_{L^1_tB_{\infty,2}^0}}.
\end{equation}
Let $N\in\NN$, then we have by Bernstein inequalities and Proposition \ref{max-pro}
\begin{eqnarray*}
\nonumber\|\theta\|_{L^1_tB_{\infty,2}^0}&\le&\|S_N\theta\|_{L^1_tB_{\infty,2}^0}+\|(\hbox{Id}-S_N)\theta\|_{L^1_t B_{\infty,1}^0}\\
\nonumber&\lesssim& 
t\|\theta^0\|_{L^\infty}\sqrt{N}+\sum_{q\geq N }\|\Delta_q\theta\|_{L^1_tL^\infty}\\
&\lesssim& \sqrt{N}\|\theta^0\|_{L^\infty} t+\sum_{q\geq N }2^{q\frac2p}\|\Delta_q\theta\|_{L^1_tL^p}.
\end{eqnarray*}
Using Theorem \ref{thm99} and $p>2$ we obtain 
\begin{eqnarray*}
\nonumber
\sum_{q\geq N-1 }2^{q\frac2p}\|\Delta_q\theta\|_{L^1_tL^p}&\lesssim&\sum_{q\geq N-1 }2^{q(\frac2p-1)}\Big(  \|\theta^0\|_{L^p}+\|\theta^0\|_{L^\infty}\int_0^t\|\omega(\tau)\|_{L^p}d\tau\Big)
\\
&\lesssim&  \|\theta^0\|_{L^p}+2^{N(-1+\frac2p)}\|\theta^0\|_{L^\infty}\int_0^t\|\omega(\tau)\|_{L^p}d\tau.
\end{eqnarray*}
Thus, we get 
\begin{eqnarray}\label{es43}
\nonumber\|\theta\|_{L^1_tB_{\infty,2}^0}\lesssim\sqrt{N}\|\theta^0\|_{L^\infty} t+ \|\theta^0\|_{L^p}+2^{N(-1+\frac2p)}\|\theta^0\|_{L^\infty}\int_0^t\|\omega(\tau)\|_{L^p}d\tau.
\end{eqnarray}
We choose $N$ as follows
 $$
 {N= \Bigg[\frac{\log\Big(e+\int_0^t\|\omega(\tau)\|_{L^p}d\tau\Big)}{(1-2/p)\log2}\Bigg]+1.}
 $$ 
 Then it follows
\begin{equation*}
\|\theta\|_{L^1_tB_{\infty,2}^0}\lesssim\|\theta^0\|_{L^\infty\cap L^p}+\|\theta^0\|_{L^\infty} t\log^{\frac12}\Big(e+\int_0^t\|\omega(\tau)\|_{L^p}d\tau\Big).
\end{equation*}
Combining this estimate with (\ref{gr1}) we get
\begin{eqnarray*}
\|\theta\|_{L^1_tB_{\infty,2}^0}&\lesssim& \|\theta^0\|_{L^\infty\cap L^p}+\|\theta^0\|_{L^\infty} t\, \log^{\frac12}\Big(e+C_0e^{C_0 t}e^{C\|\theta\|_{L^1_tB_{\infty,2}^0}}\Big)\\
&\le&C_0(1+t^2)+C\|\theta^0\|_{L^\infty} t\|\theta\|_{L^1_tB_{\infty,2}^0}^{\frac12}.
\end{eqnarray*}

Thus we get for every $t\in\RR+$
\begin{equation*}
\|\theta\|_{L^1_tB_{\infty,2}^0}\le C_0 (1+t^2).
\end{equation*}

It follows from (\ref{gr1})
\begin{equation}
\label{gr6}
\|\omega(t)\|_{L^p}\le \Phi_1(t).
\end{equation}
Applying Theorem \ref{thm99} and (\ref{gr6}) we  get 
\begin{equation}
\label{gr4}
2^q\|\Delta_q\theta\|_{L^1_tL^p}\leq \Phi_1(t),\qquad  \forall q\in\NN
\end{equation}
and thus
$$ \|\theta \|_{\widetilde{L}^1_{t} B^1_{{p, \infty}} } \leq \Phi_{1}(t).$$
This ends the proof of Proposition \ref{max-pro}.
\begin{rema} It is not hard to see  that from \eqref{gr4} one can obtain that  for every $s<1$
\begin{equation}
\label{x5}
\|\theta\|_{ L^1_tB_{p,1}^s}\leq  \|\theta\|_{ \widetilde{ L}^1_tB_{p,\infty}^1} \leq  \Phi_1(t).
\end{equation}
Combined with  Bernstein inequalities and the fact that $p>2$ this yields
\begin{equation}
\label{x9}
\|\theta\|_{ L^1_tB_{\infty,1}^\epsilon}\leq
 \Phi_1(t),
\end{equation}
for every $\epsilon <1-\frac2p$.
\end{rema}
\vspace{0.5cm}
\end{proof}

We aim now  at giving an $L^\infty$-bound of the vorticity.
 
\begin{prop}
\label{pr0} 
Let $(v,\theta)$ be a smooth solution of  \eqref{Bouss}  such that   $\omega^0,\theta^0\in L^p\cap L^\infty$ and $\mathcal{R}\theta^0\in L^\infty$, with $2<p<\infty.$  Then
we have
\begin{equation}
\label{x10}
\|\omega(t)\|_{L^\infty}+\|\mathcal{R}\theta(t)\|_{L^\infty}\le \Phi_2(t)
\end{equation}
and 
\begin{equation}
\label{x11}
\|v(t)\|_{L^\infty}\le \Phi_3(t).
\end{equation}
\end{prop}
\begin{proof}

\

$\bullet $ {\sl Proof of \eqref{x10}.}
By using the  maximum principle for the transport equation \eqref{bm1},  we get
$$
\|\Gamma(t)\|_{L^\infty}\leq\|\Gamma^0\|_{L^\infty}+\int_0^t\|[\mathcal{R}, v\cdot\nabla]\theta(\tau)\|_{L^\infty}
d\tau.
$$
Since the function $\mathcal{R}\theta$ satisfies the equation
\begin{equation}
\label{rtheta}
\big(\partial_t+v\cdot\nabla+|\DD|\big)\mathcal{R}\theta=-[\mathcal{R},v\cdot\nabla]\theta,
\end{equation}
we  get by using Proposition \ref{propmax} for $p=\infty$ that 
$$
\|\mathcal{R}\theta(t)\|_{L^\infty}\leq\|\mathcal{R}\theta(t)\|_{L^\infty}+\int_0^t\|[\mathcal{R}, v\cdot\nabla]\theta(\tau)\|_{L^\infty}
d\tau.
$$
Combining the last  two estimates yields
\begin{eqnarray*}
\|\Gamma(t)\|_{L^\infty}+\|\mathcal{R}\theta(t)\|_{L^\infty}&\leq&\|\Gamma^0\|_{L^\infty}+\|\mathcal{R}\theta^0\|_{L^\infty}+2\int_0^t\|[\mathcal{R}, v\cdot\nabla]\theta(\tau)\|_{L^\infty}
d\tau
\\
&\leq & C_0+\int_0^t  \|[\mathcal{R}, v\cdot\nabla]\theta(\tau)\|_{B_{\infty,1}^0}
d\tau.
\end{eqnarray*}
It follows from  the second estimate of Theorem \ref{propcom}  and Proposition \ref{max-pro}   
\begin{eqnarray*}
\|\omega(t)\|_{L^\infty}+\|\mathcal{R}\theta(t)\|_{L^\infty}&\lesssim&C_0+\int_0^t\|\omega(\tau)\|_{L^\infty\cap L^p}\big(\|\theta(\tau)\|_{B_{\infty,1}^\epsilon}+\|\theta(\tau)\|_{L^p}\big)
d\tau
\\
&\lesssim&C_0+\|\omega\|_{L^\infty_t L^p}\big(\|\theta\|_{L^1_tB_{\infty,1}^\epsilon}+t\|\theta^0\|_{L^p}\big)\\
&+&\int_0^t\|\omega(\tau)\|_{L^\infty}\big(\|\theta(\tau)\|_{B_{\infty,1}^\epsilon}+\|\theta^0\|_{L^p}\big)
d\tau.
\end{eqnarray*}
Let $0<\epsilon<1-\frac2p$ then using \eqref{x9} we get
$$
\|\omega(t)\|_{L^\infty}+\|\mathcal{R}\theta(t)\|_{L^\infty}\lesssim\Phi_1(t)+\int_0^t\|\omega(\tau)\|_{L^\infty}\big(\|\theta(\tau)\|_{B_{\infty,1}^\epsilon}+\|\theta^0\|_{L^p}\big)d\tau.
$$
 Therefore we obtain by the   Gronwall lemma  and a new use of \eqref{x9} that 
 \begin{eqnarray*}
\|\omega(t)\|_{L^\infty}+\|\mathcal{R}\theta(t)\|_{L^\infty}\leq \Phi_2(t).
\end{eqnarray*}

\

$\bullet $ {\sl Proof of \eqref{x11}.}

\

Let  $N\in\NN$ to be chosen later.  Using the fact that $\|\dot\Delta_qv\|_{L^\infty}\approx 2^{-q}\|\dot\Delta_q\omega\|_{L^\infty}$, we then have 
\begin{eqnarray*}
\|v(t)\|_{L^\infty}&\leq& \|\chi(2^{-N}|\DD|)v(t)\|_{L^\infty}+\sum_{q\geq-N}2^{-q}\|\dot\Delta_q\omega(t)\|_{L^\infty}\\
\\
&\le&  \|\chi(2^{-N}|\DD|)v(t)\|_{L^\infty}+2^N\|\omega(t)\|_{L^\infty}.
\end{eqnarray*} 
Applying the frequency localizing operator to the velocity equation we get
$$
\chi(2^{-N}|\DD|)v=\chi(2^{-N}|\DD|)v_0+\int_0^t\mathcal P\chi(2^{-N}|\DD|)\theta(\tau)d\tau+\int_0^t\mathcal P\chi(2^{-N}|\DD|){\rm div}(v\otimes v)(\tau)d\tau.
$$
where $\mathcal P$ stands for Leray projector. From   Bernstein inequalities, Calder\'on-Zygmund estimate and the uniform boundness of $\chi(2^{-N}|\DD|)$  we get 
\begin{eqnarray*}
\int_0^t\|\chi(2^{-N}|\DD|)\mathcal P\theta(\tau)\|_{L^\infty}d\tau&\lesssim&2^{-N\frac2p} \int_0^t\|\theta(\tau)\|_{L^p}d\tau\\
&\lesssim&t\|\theta^0\|_{L^p}.
\end{eqnarray*}
Using Proposition \ref{cor1}-(2) we find
$$
\int_0^t\|\mathcal P\chi(2^{-N}|\DD|){\rm div}(v\otimes v)(\tau)\|_{L^\infty}d\tau\lesssim 2^{-N}\int_0^t\|v(\tau)\|_{L^\infty}^2d\tau.
$$
The outcome is 
\begin{eqnarray*}
\|v(t)\|_{L^\infty}&\lesssim& \|v^0\|_{L^\infty}+t \| \theta_0\|_{L^p}+ 2^{-N}\int_0^t \| v(\tau)\|_{L^\infty}^2d\tau+2^N\|\omega(t)\|_{L^\infty}
\\
&\lesssim&   2^{-N}\int_0^t \| v(\tau)\|_{L^\infty}^2d\tau+  2^N  \Phi_2(t) 
\end{eqnarray*} 
Choosing judiciously $N$ we find
$$
\|v(t)\|_{L^\infty}\le\Phi_2(t)\Big(1+\Big(\int_0^t\|v(\tau)\|_{L^\infty}^2d\tau\Big)^{\frac12}\Big).
$$From Gronwall lemma we get
$$
\|v(t)\|_{L^\infty}\le \Phi_3(t).
$$

\

\end{proof}

Now we will describe the last part of the {\it a priori} estimates.  Following the program exposed in the beginning,
the aim is to get estimates on  $\nabla v$. 
\begin{prop}
\label{pr10}
Let $(v,\theta)$ be a smooth solution of \eqref{Bouss} and  $\omega^0,\theta^0\in B_{\infty,1}^0\cap L^p$  with $p\in]2,\infty[.$ Then
$$
\|\theta(t)\|_{B_{\infty,1}^0}+\|\omega(t)\|_{B_{\infty,1}^0}+\|v(t)\|_{B_{\infty,1}^1}\leq \Phi_3(t).
$$

\end{prop}
\begin{proof} 
By using  the logarithmic estimates of 
Theorem \ref{thmlog} for the equations \eqref{bm} and \eqref{rtheta},  we obtain
\begin{equation}\label{eqlog}
\|\Gamma(t)\|_{B_{\infty,1}^0}+\|\mathcal R\theta(t)\|_{B_{\infty,1}^0} \lesssim\Big(C_0+\big\|[\mathcal{R}, v\cdot\nabla]\theta\big\|_{L^1_tB_{\infty,1}^0}\Big)\Big(1+\|\nabla v\|_{L^1_tL^\infty}  \Big).
\end{equation}
Thanks to  Theorem \ref{propcom}, Propositions \ref{pr0}, \ref{max-pro} and \eqref{x9} we get
\begin{eqnarray*}
\big\|[\mathcal{R}, v\cdot\nabla]\theta\big\|_{L^1_tB_{\infty,1}^0}&\lesssim &\int_0^t(\|\omega(\tau)\|_{L^\infty}+\|\omega(\tau)\|_{L^p})\big(\|\theta(\tau)\|_{B_{\infty,1}^\epsilon}+\|\theta(\tau)\|_{L^p}\big)d\tau
\\
&\lesssim & \Phi_2(t).
\end{eqnarray*}
By easy computations  we get
\begin{eqnarray}
\nonumber
\|\nabla v\|_{L^\infty}&\le& \|\nabla\Delta_{-1}v\|_{L^\infty}+\sum_{q\in\NN}\|\Delta_q\nabla v\|_{L^\infty}
\\
\nonumber
&\lesssim& \|\omega\|_{L^p}+\sum_{q\in\NN}\|\Delta_q\omega\|_{L^\infty}
\\
\label{x8}
&\lesssim&  \Phi_1(t)+\|\omega(t)\|_{B_{\infty,1}^0}.
\end{eqnarray}
Putting together \eqref{eqlog} and \eqref{x8} leads to
\begin{eqnarray*}
\|\omega(t)\|_{B_{\infty,1}^0}\leq \|\Gamma(t)\|_{B_{\infty,1}^0}+\|\mathcal R\theta(t)\|_{B_{\infty,1}^0}\leq  \Phi_2(t)\Big(1+\int_0^t\|\omega(\tau)\|_{B_{\infty,1}^0}d\tau\Big).
\end{eqnarray*}
Thus we obtain from Gronwall inequality
\begin{equation}\label{x12}
\|\omega(t)\|_{B_{\infty,1}^0}\leq\Phi_3(t).
\end{equation}
Coming back to \eqref{x8} we get
\begin{eqnarray*}
\nonumber
\|\nabla v(t)\|_{L^\infty}\leq   \Phi_3(t).
\end{eqnarray*}
Let us move to the estimate of $v$ in the space $B_{\infty,1}^1$. By definition we have
$$
\|v(t)\|_{B_{\infty,1}^1}\lesssim \|v(t)\|_{L^\infty}+\|\omega(t)\|_{B_{\infty,1}^0}.
$$
Combined with \eqref{x11} and   \eqref{x12} this yields 
$$
\|v(t)\|_{B_{\infty,1}^1}\leq \Phi_3(t)
$$
as claimed. To estimate $\|\theta(t)\|_{B_{\infty,1}^0}$ we use Theorem \ref{thmlog} and the Lipschitz estimate of the velocity.
\begin{eqnarray*}
\|\theta(t)\|_{B_{\infty,1}^0}&\lesssim& \|\theta^0\|_{B_{\infty,1}^0}\Big( 1+\int_0^t\|\nabla v(\tau)\|_{L^\infty}d\tau\Big)\\
&\lesssim& \Phi_3(t).
\end{eqnarray*}
The proof of \mbox{Proposition \ref{pr10}} is now achieved.  
\end{proof}
\begin{rema}
\label{rema7}Notice that the {\it a priori} estimates above  imply that  $\nabla\theta\in L^{1}_{\rm loc}(\RR_+;L^\infty)$. Indeed, we can   establish the following estimate by combining the smoothing effect of the temperature equation with the logarithmic estimate described by Proposition \ref{thmlog},
$$
\|\theta\|_{L^1_tB_{\infty,1}^1}\le C\|\theta^0\|_{B_{\infty,1}^0}(1+\|\nabla v\|_{L^1_tL^\infty}^2).
$$
Although we expect to use this estimate for the uniqueness part, it seems that this is not sufficient for our purpose and some technical problems arise due to the fact that Riesz transforms  do not  map continuously $L^\infty$ to itself. The crucial information  that we need for the uniqueness is  $\theta\in\widetilde L_{\textnormal{loc}}^1(\RR_+; {B}_{p, \infty}^{1}).$ 
\end{rema}

\subsection{Uniqueness}
We now  show that the Boussinesq system \eqref{Bouss} has a unique solution in the class 
$$\mathcal{X}_T= (L^\infty_TB_{\infty,1}^0\cap L^1_TB_{\infty,1}^1)\times (L^\infty_TL^p\cap \widetilde L^1_T B_{p,\infty}^{1}),\quad 2<p<\infty.
$$
Let $(v^1,\theta^1)$ and  $(v^2,\theta^2)$  two solutions of \eqref{Bouss} belonging to the space $\mathcal{X}_T,$ and denote 
$$v=v^2-v^1,\quad\theta=\theta^2-\theta^1.
$$ Then we have the equations
\begin{equation*}
\left\{ 
\begin{array}{ll}
 \partial_t v+v^2\cdot\nabla v=-\nabla p-v\cdot\nabla v^1+\theta e_2\\ 
\partial_t\theta+v^2\cdot\nabla \theta+|\textnormal{D}|\theta=-v\cdot\nabla\theta^1\\
v_{| t=0}=v^0, \quad \theta_{| t=0}=\theta^0. 
\end{array} \right.
\end{equation*}
According to  Theorem \ref{thmlog} we have
$$
\|v(t)\|_{B_{\infty,1}^0}\lesssim (1+V_1(t))\Big(\|v^0\|_{B_{\infty,1}^0}+\|\nabla p\|_{L^1_tB_{\infty,1}^0}+\| v\cdot\nabla v^1\|_{L^1_tB_{\infty,1}^0}+\|\theta\|_{L^1_tB_{\infty,1}^0}\Big),
$$
with $V_1(t)=\|\nabla v^1\|_{L^1_tL^\infty}.$
A straightforward calculus using the incompressibility of the flows gives
\begin{eqnarray*}
\nabla p&=&-\nabla\Delta^{-1}\textnormal{div }(v\cdot\nabla(v^1+v^2))+\nabla\Delta^{-1}\partial_2\theta\\
&=&\hbox{I}+\hbox{II}.
\end{eqnarray*}
To estimate the first term of the RHS we use the definition
$$
\|\hbox{I}\|_{B_{\infty,1}^0}\lesssim\|(\nabla\Delta^{-1}\textnormal{div })\textnormal{div }\Delta_{-1}(v\otimes (v^1+v^2))\|_{L^\infty}+\|v\cdot\nabla(v^1+v^2)\|_{B_{\infty,1}^1}
$$
From Proposition \ref{cor1}-(2) and Besov embeddings  we have
\begin{eqnarray*}
\|(\nabla\Delta^{-1}\textnormal{div })\textnormal{div }\Delta_{-1}(v\otimes (v^1+v^2))\|_{L^\infty}&\lesssim& \|v\otimes (v^1+v^2)\|_{L^\infty}\\
&\lesssim&\|v\|_{B_{\infty,1}^0} \|v^1+v^2\|_{B_{\infty,1}^0}.
\end{eqnarray*}
Using the incompressibility of $v$ and using Bony's decomposition one can easily obtain
$$
\|v\cdot\nabla(v^1+v^2)\|_{B_{\infty,1}^0}\lesssim \|v\|_{B_{\infty,1}^0}\|v^1+v^2\|_{B_{\infty,1}^1}.
$$
Putting together  these estimates yields
\begin{equation}\label{s12}
\|\hbox{I}\|_{B_{\infty,1}^0}\lesssim \|v\|_{B_{\infty,1}^0}\|v^1+v^2\|_{B_{\infty,1}^1}.
\end{equation}
Let us now show how to estimate the second term $\hbox{II}$. By using Besov embeddings and 
 Calder\'on-Zygmund estimate we get
 \begin{eqnarray*}
 \|\hbox{II}\|_{B_{\infty,1}^0}&\lesssim& \|\nabla\Delta^{-1}\partial_2\theta\|_{B_{p,1}^{\frac2p}}\\
 &\lesssim&  \|\theta\|_{B_{p,1}^{\frac2p}}.
 \end{eqnarray*}
 Combining this estimate with (\ref{s12}) yields
 \begin{eqnarray}
 \label{1170}
\nonumber\|v(t)\|_{B_{\infty,1}^0}\lesssim e^{CV(t)}\Big(\|v^0\|_{B_{\infty,1}^0}&+&\int_0^t\| v(\tau)\|_{B_{\infty,1}^0}\big[1+\|(v^1,v^2)(\tau)\|_{B_{\infty,1}^1}\big]d\tau\Big)\\&+&e^{CV(t)}
\|\theta\|_{L^1_tB_{p,1}^{\frac2p}},
 \end{eqnarray}
where $V(t):=\|(v^1,v^2)\|_{L^1_tB_{\infty,1}^1}.$

We have now to estimate $\|\theta\|_{L^1_tB_{p,1}^{\frac2p}}$.  By applying $\Delta_q$ to the equation of $\theta $ and arguing similarly to the proof of Theorem \ref{thm99} we obtain
\begin{eqnarray*}
\nonumber\|\theta_q(t)\|_{L^p} \lesssim   e^{-ct2^q}\|\theta_q^0\|_{L^p}&+&\int_0^t e^{-c2^q(t-\tau)} \|\Delta_q(v\cdot\nabla \theta^1)(\tau)\|_{L^p}d\tau\\
&+&\int_0^te^{-c2^q(t-\tau)} \|\big[v^2\cdot\nabla,\Delta_q\big]\theta(\tau) \|_{L^p}d\tau.
\end{eqnarray*}
Remark, first, that an obvious H\"older inequality yields that for every $\varepsilon\in[0,1]$ there exists an absolute constant $C$ such that 
$$
\int_0^t e^{-c\tau2^q}d\tau\leq C t^{\varepsilon}2^{-q(1-\varepsilon)},\qquad \forall\, t\geq 0.
$$
Using this fact  and integrating in time 
\begin{eqnarray}
\label{113}
\nonumber2^{q\frac2p}\|\theta_q\|_{L^1_tL^p}& \lesssim &  2^{q(-1+\frac2p)}\|\theta_q^0\|_{L^p}\\
\nonumber&+& t^{\varepsilon}2^{q(-1+\varepsilon+\frac2p)}\int_0^t \Big( \|\Delta_q(v\cdot\nabla \theta^1)(\tau)\|_{L^p}+ \|[v^2\cdot\nabla,\Delta_q]\theta(\tau) \|_{L^p}\Big)d\tau\\
&=& 2^{q(-1+\frac2p)}\|\theta_q^0\|_{L^p}+{\rm I}_q(t)+{\rm II}_q(t).
\end{eqnarray}
Using Bony's decomposition we get easily
\begin{eqnarray*}
 \|\Delta_q(v\cdot\nabla \theta^1)(t)\|_{L^p}& \lesssim&   \|v(t)\|_{L^\infty}\sum_{j\leq q+2} 2^{j}  \|\Delta_j\theta^1(t)\|_{L^p}\\
 &+& 2^q \|v(t)\|_{L^\infty}\sum_{j\geq q-4}   \|\Delta_j\theta^1(t)\|_{L^p}.
  \end{eqnarray*}
  Integrating in time we get
  \begin{eqnarray}\label{114}
 \nonumber {\rm{I}}_q (t) & \lesssim& t^{\varepsilon} \|v\|_{L^\infty_tL^\infty}2^{q(-1+\frac2p+\varepsilon)}(|q|+1)\|\theta^1\|_{\widetilde L^1_tB_{p,\infty}^1}\\ 
 \nonumber&+& t^{\varepsilon} \|v\|_{L^\infty_tL^\infty}2^{q(\frac2p+\varepsilon)}\sum_{j\geq q-4}   \|\Delta_j\theta^1\|_{L^1_tL^p}\\
 &\lesssim&t^{\varepsilon} \|v\|_{L^\infty_tL^\infty}2^{q(-1+\frac2p+\varepsilon)}(|q|+1)\|\theta^1\|_{\widetilde L^1_tB_{p,\infty}^1}.
  \end{eqnarray}
  To estimate the term ${\rm II}_q$ we use the following classical commutator $(2/p<1)$
  $$
  \|[v^2\cdot\nabla,\Delta_q]\theta \|_{L^p}\lesssim  2^{-q\frac2p}\|\nabla v^2\|_{L^\infty}\|\theta\|_{B_{p,1}^{\frac2p}}.
  $$
  
  Thus we obtain,
\begin{equation}\label{115}
{\rm II}_q (t) \lesssim t^{\varepsilon}2^{q(-1+\varepsilon)}\|\nabla v^2\|_{L^\infty_tL^\infty}\|\theta\|_{L^1_tB_{p,1}^{\frac2p}}
\end{equation}

We choose $\varepsilon>0$ such that $-1+\frac2p+\varepsilon<0,$ which is possible since $p>2.$ Then  summing (\ref{113}) and using (\ref{114}) and (\ref{115}) we get 
\begin{eqnarray*}
\|\theta\|_{L^1_tB_{p,1}^{\frac2p}}&\lesssim & \|\theta^0\|_{L^p}+t^{\varepsilon} \|v\|_{L^\infty_tL^\infty}\|\theta^1\|_{\widetilde L^1_tB_{p,\infty}^1}+t^{\varepsilon}\|\nabla v^2\|_{L^\infty_tL^\infty}\|\theta\|_{L^1_tB_{p,1}^{\frac2p}}.
\end{eqnarray*}
For small $t\in[0,\delta]$ one can obtain
$$
\|\theta\|_{L^1_tB_{p,1}^{\frac2p}}\lesssim  \|\theta^0\|_{L^p}+t^{\varepsilon} \|v\|_{L^\infty_tL^\infty}\|\theta^1\|_{\widetilde L^1_tB_{p,\infty}^1}.
$$ 
Plugging this estimate into \eqref{1170} we find
$$
\|v\|_{L^\infty_tB_{\infty,1}^0}\lesssim e^{CV(t)}\Big(\|v^0\|_{B_{\infty,1}^0}+\|\theta^0\|_{L^p}+t\| v\|_{L^\infty_tB_{\infty,1}^0}+ t^{\varepsilon} \|v\|_{L^\infty_tL^\infty}\|\theta^1\|_{\widetilde L^1_tB_{p,\infty}^1}\Big).
$$
If $\delta$ is sufficiently small then we get for $t\in[0,\delta]$
\begin{equation}\label{uni1}
\|v\|_{L^\infty_tB_{\infty,1}^0}\lesssim \|v^0\|_{B_{\infty,1}^0}+\|\theta^0\|_{L^p}.
\end{equation}
This gives in turn
\begin{equation}\label{uni2}
\|\theta\|_{L^1_tB_{p,1}^{\frac2p}}\lesssim \|v^0\|_{B_{\infty,1}^0}+\|\theta^0\|_{L^p}.
\end{equation}
This gives in particular the uniqueness on  $[0,\delta]$. An iteration of this argument  yields   the uniqueness in $[0,T]$. 

\subsection{Existence}
We consider the following system
\begin{equation} 
\left\{ \begin{array}{ll} 
\partial_{t}v_n+v_n\cdot\nabla v_n+\nabla p_n=\theta_n e_{2}\\ 
\partial_{t}\theta_n+v_n\cdot\nabla\theta_n+\vert \textnormal{D}\vert\theta_n=0\\
\textnormal{div}v_n=0\\
{v_n}_{| t=0}=S_nv^{0}, \quad {\theta_n}_{| t=0}=S_n\theta^{0}  
\end{array} \right. \tag{B$_{n}$}
\end{equation}
By using the same method as \cite{hk} we can prove that this system has a unique  local smooth solution $(v_n,\theta_n)$. The global existence of these solutions is governed by  the following criterion: we can push the construction beyond the time $T$ if  the quantity $\|\nabla v_n\|_{L^1_T L^\infty}$ is finite. Now from the {\it a priori} estimates the Lipschitz norm can not blow up in finite time and then the solution $(v_n,\theta_n)$ is globally defined. Once again from the {\it a priori} estimates   we have  for $2<p<\infty$ 
$$
\|v_n\|_{ L^\infty_TB_{\infty,1}^1}+\|\omega_n\|_{L^\infty_T L^p}+\|\theta_n\|_{L^\infty_T(B_{\infty,1}^{0}\cap L^p)}\leq \Phi_3(T).
$$
Then it follows that up to an extraction the sequence $(v_n,\theta_n)$ is weakly convergent to  $(v,\theta)$ belonging to $L^\infty_TB_{\infty,1}^1\times  L^\infty_T(B_{\infty,1}^{0}\cap L^p),$ with $\omega\in L^\infty_T L^p$.  For $(n,m)\in\NN^2$ we set $v_{n,m}=v_n-v_m$ and $\theta_{n,m}=\theta_n-\theta_m$ then  according to the estimate (\ref{uni1}) and (\ref{uni2}) we get  for $T=\delta$
$$\|v_{n,m}\|_{L^\infty_TB_{\infty,1}^0}+\|\theta_{n,m}\|_{L^1_TB_{p,1}^{\frac2p}}\leq \Phi_3(t) \big( \|S_nv^0-S_m v^0\|_{B_{\infty,1}^0}+\|S_n\theta^0-S_m \theta^0\|_{L^p}\big).
$$
This shows that $(v_n,\theta_n)$ is a Cauchy sequence in the Banach space $L^\infty_TB_{\infty,1}^0\times L^1_TB_{p,1}^{\frac2p}$ and then it converges strongly to $(v,\theta).$ This allows to pass to the limit in the system $({\rm B}_n)$ and then we get that $(v,\theta)$ is a solution of the Boussinesq system (\ref{Bouss}).

\section*{Appendix: proof of  Lemma \ref{propDelta}}
\label{technical}

 We have from Bony's decomposition
\begin{eqnarray*}
[\Delta_q, v\cdot\nabla]\theta&=&\sum_{|j-q|\le 4}[\Delta_q, S_{j-1}v\cdot\nabla]
\Delta_j\theta+\sum_{|j-q|\le4}[\Delta_q, \Delta_jv\cdot\nabla]S_{j-1}\theta\\
\nonumber&+&\sum_{j\geq q-4} [\Delta_q, 
\Delta_jv\cdot\nabla]\widetilde{\Delta}_j\theta\\
&:=& \mbox{I}_q+\mbox{II}_q+\mbox{III}_q.
\end{eqnarray*}
Observe first that 
 $$
\mbox{I}_q=\sum_{|j-q|\le 4}h_q\star( S_{j-1}v\cdot\nabla\Delta_j\theta)-S_{j-1}v\cdot(h_q\star \nabla\Delta_j\theta)
$$
where  $\hat h_q(\xi)=\varphi({2^{-q}}\xi)$. Thus, Lemma \ref{commu} and Bernstein inequalities yield
\begin{eqnarray*}
\nonumber\|\mbox{I}_q\|_{L^p}&\lesssim&\sum_{|j-q|\le4}\|xh_q\|_{L^1}\|\nabla S_{j-1} 
v\|_{L^p}\|\nabla\Delta_j\theta\|_{L^\infty}\\
\nonumber
&\lesssim&\|\nabla  v\|_{L^p}\|h_0\|_{L^1}\sum_{|j-q|\le4}2^{j-q}\|\Delta_j\theta\|_{L^\infty}\\
&\lesssim&\|\nabla v\|_{L^p}\|\theta\|_{B_{\infty,\infty}^0}.
\end{eqnarray*}
To estimate  the second term we use once again Lemma \ref{commu}
\begin{eqnarray*}
\|\hbox{II}_q\|_{L^p}&\lesssim&\sum_{|j-q|\le 4}2^{-q}\|\Delta_j\nabla v\|_{L^p}\|\nabla S_{j-1}\theta\|_{L^\infty}\\&\lesssim&
 \|\nabla v\|_{L^p}\sum_{|j-q|\le 4\atop k\le j-2}2^{k-q}\|\Delta_k\theta\|_{L^\infty}\\
&\lesssim&\|\nabla v\|_{L^p}\|\theta\|_{B_{\infty,\infty}^0}.
\end{eqnarray*}
 Let us now move to the remainder term. We separate it into two terms: high frequencies and low frequencies.
\begin{eqnarray*}
\hbox{III}_q&=&\sum_{j\geq q-4\atop j\in\NN} [\Delta_q\partial_i, \Delta_jv^i]\widetilde{\Delta}_j\theta+ [\Delta_q, \Delta_{-1}v\cdot\nabla]\widetilde{\Delta}_{-1}\theta\\
&:=&\hbox{III}_q^1+\hbox{III}_q^2.
\end{eqnarray*}
For the first term we don't need to use the structure of the commutator. We estimate separately each term of the commutator by using Bernstein inequalities.
\begin{eqnarray*}
\|\hbox{III}_q^1\|_{L^p}&\lesssim& \sum_{j\geq q-4\atop j\in\NN} 2^q\|\Delta_jv\|_{L^p}\|\widetilde\Delta_j\theta\|_{L^\infty}\\
&\lesssim&
 \|\nabla v\|_{L^p}\sum_{j\geq q- 4}2^{q-j}\|\widetilde\Delta_j\theta\|_{L^\infty}\\
&\lesssim&\|\nabla v\|_{L^p}\|\theta\|_{B_{\infty,\infty}^0}.
\end{eqnarray*}
For the second term we use Lemma \ref{commu} combined with Bernstein inequalities.
\begin{eqnarray*}
\|\hbox{III}_q^2\|_{L^p}&\lesssim&\|\nabla\Delta_{-1} v\|_{L^p}\|\nabla\widetilde\Delta_{-1}\theta\|_{L^\infty}\\
&\lesssim&\|\nabla v\|_{L^p}\|\theta\|_{B_{\infty,\infty}^0}.
\end{eqnarray*}

\end{document}